\documentclass[reqno,11pt]{amsart} 
\usepackage[latin1]{inputenc}
\usepackage{amsmath}  
\usepackage{amsfonts} 
\usepackage[french,english]{babel}
\usepackage{amssymb} 
\usepackage{dsfont}\let\mathbb\mathds
\usepackage[all,dvips]{xy}
\usepackage{latexsym}
\usepackage[dvips]{graphicx}
\newenvironment{Dem}{\textit{Proof :}}{\begin{flushright}$\Box$\end{flushright}} 
\newenvironment{DemP}{\textit{Proof of the Proposition \ref{divisor}:}}{\begin{flushright}$\Box$\end{flushright}} 
\newenvironment{DemT}{\textit{Proof of the Theorem \ref{esencial}:}}{\begin{flushright}$\Box$\end{flushright}}

\newtheorem{The}{Theorem}[section] 
\newtheorem{Pro}{Proposition}[section] 
\newtheorem{Lem}[The] {Lemma}

\theoremstyle{definition}

\theoremstyle{remark} 
\newtheorem{Rmk}[The]{Remark} 
 
\numberwithin{equation}{section}

\newcommand{\fs}{\mathcal{O}_C}
\newcommand{\fc}{\mathcal{W}_C}
\newcommand{\trans}{\mathcal{O}_{P_0}(t^*_{-y}\Xi_0) }
\newcommand{\Trans}{\mathcal{O}_{P_0}(t^*_y\Xi_0) }

\newcommand{\im}{\mathop{\rm Im}\nolimits}

\newcommand{\pic}{\mathop{\rm Pic}\nolimits}
\newcommand{\supp}{\mathop{\rm Supp}\nolimits}
\newcommand{\Div}{\mathop{\rm Div}\nolimits}

\newcommand{\Ker}{\mathop{\rm Ker}\nolimits}
\newcommand{\End}{\mathop{\rm End}\nolimits}
\newcommand{\Deg}{\mathop{\rm deg}\nolimits}
\newcommand{\Dim}{\mathop{\rm dim}\nolimits}

\begin{document} 
\title{Prym--Tyurin Varieties coming from correspondences with fixed points.} 
 
\author{Angela Ortega} 
\address{Angela Ortega, Instituto de Matem\'aticas, U.N.A.M., Morelia , M\'exico} 
\email{ortega@matmor.unam.mx} 
 
\subjclass[2000]{}
 
\date{\today.} 
 
 
\begin{abstract} 
Our main theorem is an improvement of the Criterion of Kanev about Prym--Tyurin varieties induced by correspondences, which includes
correspondences with fixed points. We give some examples of Prym--Tyurin varieties using this criterion.  
\end{abstract} 
 
\maketitle 
 
\markboth{Angela Ortega}{Prym--Tyurin Varieties}

\section{Introduction}

Let  $C$ be a smooth projective curve over $\mathbb{C}$ and let $J=JC$ be the Jacobian of $C$ with canonical
principal polarization $\Theta$. Consider $P \stackrel{i}{\hookrightarrow} J$ an abelian subvariety with
induced polarization $i^*\Theta$. We say that $P$ is a Prym--Tyurin variety of exponent $q$ for the curve
$C$ if it exists a principal polarization $\Xi \subset P$ such that
\begin{displaymath}
i^*\Theta \equiv q\Xi.
\end{displaymath}
 In~\cite{Kanev} Kanev shows the following sufficient condition for
a subvariety of $J$ to be a Prym--Tyurin variety for $C$. If $\gamma$ is an endomorphism of $J$ induced by an effective
fixed point free symmetric correspondence which satisfies the equation 
\begin{equation} \label{equation}
(1-\gamma)(\gamma +q-1)=0, \qquad q \geq 2,
\end{equation}
then $P=\im (1-\gamma)$ is a Prym--Tyurin variety of exponent q.
However, this criterion does not include the case of a Prym variety (of exponent 2) associated to a double cover
over $C$ with two branch points (cf.~\cite{Mum}). In this work we give some conditions to extend this criterion to a correspondence 
with fixed points.
\\
A correspondence is a line bundle $L$ on $C\times C$ (alternatively a divisor on $C\times C$ ); we say that the correspondence is effective if it is 
defined by an effective divisor on $C\times C$. Throughout this paper we will consider effective correspondences. For any point $p \in C$ define the line bundle on C
\begin{displaymath}
L(p):= L_{|_{\{ p\}\times C}}.
\end{displaymath} 
A point of $p \in C$, with $L(p) = \fs (D)$, is a fixed point of $L$ if $D-p$ is an effective divisor. We prove the following 
\begin{The}\label{esencial} 
Let $L$ be an effective symmetric correspondence on $C \times C$ of bidegree $(d,d)$  with $2n$ fixed points and $\gamma_L$ the endomorphism of the 
Jacobian $J$ induced by  $L$. Suppose $n \leq d$ and
\\
a) $\gamma_L$ satisfies the equation ~\ref{equation}. 
\\
b)There are  $n$  distinct fixed points $p_1,\ldots , p_n \in C$ such that 
\begin{displaymath}
p_1, \ldots ,p_i \in D_i, \qquad  p_i \notin D_i - p_i, \qquad i=1, \ldots, n,  
\end{displaymath}
where $L(p_i)= \fs (D_i)$ with $D_i$ effective divisors.
\\
Then $P := \im (1-\gamma_L)$ is a Prym--Tyurin variety of exponent $q$ for the curve $C$. 
Moreover, there exist theta divisors $\Theta$ and $\Xi$ on $J$ and $P$ respectively such that $i^*\Theta = q\Xi$.
\end{The}

\begin{Rmk}
If $L = \mathcal{O}_{C\times C}(\mathcal{D})$, the condition $p_i \notin D_i - p_i$ for $ i=1, \ldots, n$ of the Theorem \ref{esencial} is equivalent to 
say that $\mathcal{D}$ intersects  transversally to the diagonal $\Delta$ in $C\times C$ at the points $p_1, \ldots , p_n$. 
\end{Rmk}
En particular, when $d=n=1$ this criterion recovers the example of Prym varieties associated to double covers with two ramification points. 
We will use this criterion to construct new examples of Prym--Tyurin varieties. 

\emph{Acknowledgments.} I am grateful to Sev\'{\i}n Recillas for his advice and encouragement, this article is dedicated to him. I would like to thank to A. S\'anchez for very useful discussions.

\section{Subvarieties of the Jacobian}

Let us recall some generalities about correspondences and abelian varieties.
Let $g$ be the genus of $C$. We will consider correspondences on $C$ as line bundles on $C \times C$ which are defined
by effective divisors on $C \times C$. The bidegree $(d_1,d_2)$ of a correspondence $L \in \pic (C \times C)$ is defined by $d_1 = 
\Deg L_{|_{ C \times \{ s \} }}$ and $d_2 = \Deg L_{|_{\{ t \}\times C}}$, and it is independent of the points $s,t \in C$.
Two correspondences $L$ and $L'$ are equivalents if there are line bundles $L_1$ and $L_2$ on $C$ such that 
\begin{displaymath}
L'=L \otimes \pi_1^* L_1 \otimes \pi_2^*L_2,
\end{displaymath} 
where $\pi_1$ and $\pi_2$ are the canonical projections of $C \times C$. A correspondence $L$ 
induces an endomorphism of $\pic (C)$ given by
\begin{displaymath}
\gamma_L: \fs(\sum r_i p_i) \mapsto L(p_1)^{r_1} \otimes \cdots \otimes  L(p_n)^{r_n},
\end{displaymath}
with $p_i \in C$ and $r_i \in \mathbb{Z}$. 
For all  $N \in \pic (C)$ we have $\Deg \gamma_L(N)=d_2 \Deg (N)$; in particular, a correspondence induces an endomorphism of the Jacobian $\gamma_L: J 
\rightarrow J$, 
which does not depend on the class of equivalence of $L$ and every endomorphism of $J$ is obtained in this way (Theorem 11.5.1. \cite{B-L}). An 
endomorphism $\gamma_L$ associated to $L$ is symmetric if and only if  $\tau^* L = L$ where $\tau: (x,y) \mapsto (y,x)$ on $C \times C$, and then 
$d_1=d_2=d$.
Consider a symmetric endomorphism $\gamma \in \End (J)$ such that $1-\gamma$ is a primitive endomorphism satisfying 
\begin{equation*} 
(1-\gamma)^2=q(1-\gamma),
\end{equation*}
for some positive integer $q$. Define 
\begin{displaymath}
P:= \im(1-\gamma),
\end{displaymath}
which is an abelian subvariety of $J$ of exponent $q$. 
Let $P \stackrel{i}{\hookrightarrow} JC$ be the inclusion and 
$i^*\Theta$ the induced polarization on $P$.
\begin{Pro} (Prop. 1.6~\cite{Kanev}, Prop.12.1.8~\cite{B-L}) \label{diag}
There exist a  principally polarized abelian variety $(P_0, \Xi_0)$,  homomorphisms 
$\sigma, j$ and an isogeny $ \mu$ in the diagram
\begin{eqnarray}
\shorthandoff{;:!?}
\xymatrix @!0 @R=1.5cm @C=1.8cm  {
        P_0  \ar[d]_{\mu} \ar@<2pt>[dr]^{j}  & \\
	P \ar@{^{(}->}[r]^{i}  & J  \ar@<2pt>[ul]^{\sigma}   
     }
\end{eqnarray}
such that they verify the following :

$a) j^*\Theta \equiv q\Xi_0$ \\
$
b) \  i \circ \mu = j ; \qquad j \circ \sigma = 1-\gamma ; \qquad \sigma \circ j= q_{P_0} ;\qquad \sigma \circ \gamma= (1-q)\sigma.
$
\end{Pro}

\section{The proof of the criterion}

The idea of the proof of the Kanev's Criterion is to find divisors $\Theta$ on $J$ and  $\Xi$ on $P$ such that we  have the equality of
divisors $i^*\Theta = q\Xi$. Our contribution to the original proof resides in the following proposition, proved in the context of a correspondence 
which fixed points verify the condition b) of the Theorem \ref{esencial}.  
\begin{Pro}\label{divisor}
There exists theta divisors $\Theta$ and $\Xi_0$ on $J$ and $P_0$ such that $j^*(\Theta)= q\Xi_0 $.
\end{Pro}
The proof of the criterion follows from this proposition.
\\
\\
\begin{DemT} 
If $\Theta$ is the divisor of the previous proposition the divisor $i^*\Theta$ is well defined, since by the Proposition \ref{diag} $j=i \circ \mu$ and 
$\mu$ is an isogeny. Put $i^*\Theta = \sum _k r_k D_k$ with $r_k > 0$ and $\mu^*D_k =\sum_l \Xi_{kl}$ with pairwise different irreducible divisors $D_k$ on 
$P$ and $\Xi_{kl}$ on $P_0$. Hence $j^*\Theta = \mu^*i^*\Theta = \sum_{k,l}r_k \Xi_{kl}$. According to the Proposition \ref{divisor} $r_k=q$ for all $k$ and
 $\Xi_0 = \sum_{k,l} \Xi_{kl}$. We define $\Xi= \sum_k D_k$. Since $\mu^*\Xi = \Xi_0$ defines a principal polarization, the isogeny $\mu$ is an isomorphism 
and $\Xi$ defines also a principal polarization. Therefore $i^*\Theta =q \Xi$.  
\end{DemT}
Let $\alpha_c : C \rightarrow J$ be the Abel--Jacobi morphism with base point $c\in C$ and $L=\mathcal{O}_{C \times C}(\mathcal{D})$ an effective 
correspondence. 
\begin{Lem} \label{epsi} 
There exists $\varepsilon \in \pic(C)$ such that for all $y \in P_0$ and for all $c \in C$ 
\begin{displaymath}
\alpha^*_c \sigma^* \Trans = j(y)^{-1} \otimes \fs(c) \otimes L(c)^{-1} \otimes \varepsilon.
\end{displaymath}
\end{Lem}
\begin{Dem}
See ~\cite{B-L} Lemma 12.9.4.
\end{Dem}
We apply the previous lemma to the correspondence $\tau^*L$ to obtain a line bundle $\eta \in \pic(C)$. 
\begin{Lem}\label{canonico}
The lines bundle  $\varepsilon$ and $\eta$ satisfy
\begin{displaymath}
\varepsilon \otimes \eta = \fc \otimes \fs(\Delta. \mathcal{D}), 
\end{displaymath}
where $\Delta.\mathcal{D} $ is the intersection of the correspondence (as divisor in $C \times C$) with the diagonal in $C\times C $ and $\fc$ is 
the canonical divisor. 
\end{Lem}

\begin{Dem}
See  \cite{Kanev3} Theorem 7.6. (p. 211).
\end{Dem}
In particular, when $L$ is symmetric we have 
\begin{displaymath}
\varepsilon^{\otimes 2} = \fc \otimes \fs(\Delta. \mathcal{D}), 
\end{displaymath}
and hence, $\Deg \varepsilon = g +n-1$ where $2n$ is the number of fixed points of $L$ with  multiplicity. 
\\
\\
We shall adapt the Kanev's proof to admit fixed points in the criterion.  Let $L \in \pic (C \times C)$ be an effective 
symmetric correspondence with fixed points $p_1, \ldots , p_{2n} \in C$. Suppose that $p_1, \ldots , p_n $ satisfy the assumptions of the Theorem 
\ref{esencial}.
\\
\\
\begin{DemP} 
Let $\beta:=\fs (p_1 + \cdots + p_n)$
and $\varepsilon' = \varepsilon \otimes \beta^{-1}$, where $\varepsilon$ is given by the Lemma \ref{epsi}.
Observe that  $\varepsilon'$ is a line bundle of degree $g-1$. Let  $\Theta = W_{g-1} - \varepsilon'$. We can assume that the divisor $\Xi_0$ of the 
Proposition \ref{diag} is symmetric. 
\\
Firstly, we shall prove that $j^{-1}(\Theta) \subset \supp{\Xi_0} $.
Let $y \in P_0- \Xi_0$. We have to show that $j(y) \notin \Theta$, i.e., $h^0(j(y)\otimes \varepsilon')=0$.
According to the lemma \ref{epsi}, we have the following expression for $ j(y) \otimes \varepsilon'$ 
\begin{displaymath}
M:= j(y) \otimes \varepsilon' = \alpha^*_c \sigma^* \trans \otimes \fs(-c) \otimes L(c) \otimes \beta^{-1}.
\end{displaymath}
Observe that $M$ does not depend on $c$. Since $\Xi_0$ is symmetric and $y \notin \Xi_0$, $c$ is not in the divisor
$\alpha^*_c \sigma^*(t^*_{-y}\Xi_0)$ for every $c \in C$. Hence, any $c$ in the open set $U= C -\{p_1, \ldots , p_{2n} \}$ is not a base point of the 
line bundle $M \otimes \fs(c) \otimes \beta \in \pic^{g+n}(C)$, since $c$ is not a base point of $L(c)$. Then
\begin{displaymath}           
h^0(M \otimes \beta)= h^0(M\otimes \beta \otimes \fs(c)) -1, 
\end{displaymath}
for all $c \in U$. We claim that  $h^0( M \otimes \beta)= n$. Suppose that $h^0( M \otimes \beta) \geq n+1$. By Riemann--Roch we have 
\begin{eqnarray*}
h^0(\fc\otimes M^{-1} \otimes \beta^{-1} \otimes \fs(-c)) & = & h^0(M \otimes \beta \otimes \fs(c) ) -n - 1 \\
&=&  h^0(M \otimes \beta ) - n.
\end{eqnarray*}
On the other hand
\begin{equation*}
h^0(\fc\otimes M^{-1} \otimes \beta^{-1})  =  h^0(M \otimes \beta ) - n, 
\end{equation*}
hence every $c \in U$ is a base point of $\fc\otimes M^{-1} \otimes \beta^{-1} $ if  $h^0( M \otimes \beta ) \geq n + 1 $, which is impossible. 
We conclude that $h^0( M \otimes \beta)=n$. 
\\
Observe that we can write 
\begin{eqnarray*}
M \otimes \beta &=& \alpha^*_{p_1} \sigma^* \trans \otimes \fs(-p_1) \otimes L(p_1) \\
&=&  \alpha^*_{p_1} \sigma^* \trans \otimes \fs(D_1- p_1).
\end{eqnarray*}
Since $p_1 \notin D_1-p_1$ and  $p_1$ is not a base point of $\alpha^*_{p_1} \sigma^* \trans$, $p_1$ is not a 
base point of $M \otimes \beta$. 
Hence, $h^0(M \otimes \beta \otimes \fs(-p_1))= n - 1$.\\
In general, we have
\begin{eqnarray*}
M \otimes \beta \otimes \fs (-p_1 - \cdots - p_i) &=& \alpha^*_{p_{i+1}} \sigma^* \trans \otimes \fs(-p_{i+1}) \otimes L(p_{i+1})\otimes\\
&& \hspace{5cm} \otimes \fs(-p_1 -\cdots - p_i) \\
&= &  \alpha^*_{p_{i+1}} \sigma^* \trans \otimes \fs(D_{i+1}-(p_1 + \cdots +  p_{i+1})),
\end{eqnarray*}
for $i = 1, \ldots , n-1$. By the assumptions on the fixed points, $ p_{i+1}$ is not a base point of $M   \otimes \beta \otimes \fs (-p_1 - \cdots - p_i) 
=M \otimes \fs (p_{i+1} + \cdots + p_n)$ and then $h^0(M \otimes \fs(p_i + \cdots + p_n ) ) = h^0(M \otimes \fs(p_{i+1} + \cdots + p_n ) ) -1$, for 
$i = 1, \ldots , n-1$ .\\
By induction on $n$, we have $h^0(M)=0.$
\\
\\
Now, since $\Xi_0 \subset P_0$ defines a principal polarization, it is of the form $\Xi_0 = \sum \Xi_k$ where the  $\Xi_k$ are irreducible divisors 
linearly independent in $NS(P_0)$.
What we have just proved implies  $j^*(\Theta) =\sum r_k \Xi_k $, with $r_k \geq 0$. By the Proposition (\ref{diag}) $j^*(\Theta) \equiv \sum q \Xi_k$, 
hence $r_k=q$ for all $k$. Therefore $j^*\Theta= q\Xi_0$ for $\Theta= W_{g-1} -\varepsilon'$ .
\end{DemP}

\section{Examples}

{\bf Coverings over hyperelliptic curves.} In this section we consider a similar construction to the  given in \cite{LRR}. Let $X$ be an hyperelliptic 
curve of genus  $g\geq 3$, with hyperelliptic involution $i$ and let $ h : X \rightarrow \mathbb{P}^1$ be the map given by the linear system $g^1_2$. 
Consider $ f: \tilde{X} \rightarrow X$ a covering of degree 3 with two ramification points, from a projective smooth irreducible curve $\tilde{X}$.
Suppose that the branch locus of $f$ does not contain ramification points of $h$. We define a new curve $C$ by using the following cartesian diagram
\begin{eqnarray} \label{cartesiano}
\shorthandoff{;:!?}
\xymatrix  @!0 @R=1.5cm @C=2.5cm  {
        C:=(f^{(2)})^{-1}(g_2^1) \ar@{^{(}->}[r] \ar[d]_{\pi= f^{(2)}_{|_C}}   &  \tilde{X}^{(2)}  \ar[d]_{f^{(2)}}\\
        \mathbb{P}^1 \simeq g^1_2   \ar@{^{(}->}[r]^j & X^{(2)}      
     }
\end{eqnarray} 
where $f^{(2)}$ is the second symmetric product of $f$ and $\pi$ is of degree 9. Let us assume $C$ smooth and irreducible. We will define the 
same correspondence on $C$ as in ~\cite{LRR}. Let $s:\tilde{X}^2 \rightarrow  \tilde{X}^{(2)}$ be the canonical map. We denote $\tilde{C} : = s^{-1}(C) \subset 
\tilde{X}^2 $. Let $p_1: \tilde{C} \rightarrow \tilde{X} $ denote the projection on the first factor. Define  
\begin{displaymath} 
\tilde{D}:= \{ (a,b) \in \tilde{C} \times \tilde{C} \mid p_1(a)= p_1(b) \}, 
\end{displaymath}
with reduced subscheme structure. This is an effective divisor on  $ \tilde{C}^2$ containing the diagonal $\tilde{\Delta}$. Define $Y:= \tilde{D} -\tilde{\Delta} $.
The divisor $D:= (s \times s)_* (Y)$ is an effective symmetric correspondence on $C$ of bidegree (4,4).\\
Set-theoretically this 
correspondence is defined as follows. Given $z \in \mathbb{P}^1$, put $h^{-1}(z)= x + ix$ and $f^{-1}(x)=\{ x_1,x_2,x_3\}$,  $f^{-1}(ix)=\{ y_1,y_2,y_3\}$. We 
denote for $i,j \in \{1,2,3\}$      
\begin{displaymath} 
P_{ij}= x_i + y_j \in C \subset \tilde{X}^{(2)}.  
\end{displaymath}
Then  $\pi^{-1} = \{ P_{ij} | i,j=1,2,3 \}$. Let 
\begin{displaymath} 
D(P_{ij})= \sum_{l=1, l\neq j}^{3} P_{il} + \sum_{k=1, k\neq i}^{3} P_{kj}.  
\end{displaymath}
This define an effective symmetric correspondence of bidegree $(4,4)$. The associated endomorphism of the Jacobian $\gamma_D$ verifies the equation 
(cf.~\cite{LRR})
\begin{displaymath} 
\gamma_D^2 + \gamma_D - 2 =0.
\end{displaymath}
By construction, $D$ is a correspondence with fixed points coming from the ramification points of $f$. More precisely, if $x \in X$ is a branch 
point of $f$ with $x_1 = x_2$ then $\pi^{-1}(h(x))= \{ 2(x_1 + y_1), 2(x_1 + y_2),  2(x_1 + y_3), x_3 + y_1, x_3 + y_2, x_3 +y_3 \}$ and  $P_{11},P_{12},P_{13}$ are the three 
fixed points on the fiber over  $h(x)$. Moreover, these points verify
\begin{eqnarray*}
P_{11} & \in & D(P_{11}),\\
P_{11}, P_{12} & \in & D(P_{12}),\\ 
P_{11}, P_{12}, P_{13} & \in & D(P_{13}),\\ 
\end{eqnarray*}
and $P_{1i} \notin D(P_{1i})-P_{1i} $ for $i=1,2,3,$ which are the conditions of the Theorem \ref{esencial}. Applying the criterion we obtain that $P:=\im(1-\gamma_D) $ 
is a Prym--Tyurin variety of exponent 3 for the curve $C$. In order to compute
his dimension we need to calculate the degree of the ramification divisor of $\pi$. On the fiber over a Weierstrass point we have 3 simple ramification 
points i.e. with index of ramification 2. If $h(x)\in 
\mathbb{P}^1$ is the image of a branch point of $f$ then the fiber $\pi^{-1}(h(x))$ contains 3 simple ramification points. Therefore, the degree of the 
ramification divisor 
of $\pi$ is $w_{\pi} = 3(2g+2) +3(2) $. Using the Riemann--Hurwitz formula we get $g_C= 3g-2$.
By the Corollary   5.3.10 and Proposition 11.5.2 (\cite{B-L})
\begin{equation} \label{dimension} 
exp(P) \dim P  = \frac{1}{2} Tr_r (1-\gamma_D)
=  (g_C -d  + \frac{1}{2}(\Delta.D)) ,
\end{equation} 
where  $exp(P)$ is the exponent of $P$ as subvariety of $JC$,  $Tr_r$ is the rational trace and $(\Delta.D)$ denotes the number of 
fixed points of the correspondence. Hence,
\begin{displaymath}
\Dim P  =  \frac{1}{3} (g_C - 4 +3 ) = g - 1 , 
\end{displaymath}
Thus we have obtained a $2g+1$-dimensional family of Prym--Tyurin varieties of dimension $g-1$ and exponent 3.
We must to show that the curve $C$ is irreducible and smooth. With a slight modification on the argument of ~\cite{LRR}, using a suitable classifying 
homomorphism of the covering $h \circ f $, is not difficult to prove the irreducibility of $C$.
\\
\begin {Lem}
The curve $C$ given by the diagram \ref{cartesiano} is smooth.
\end{Lem}
\begin{Dem}
The curve $C \subset \tilde{X}^{(2)}$ is smooth in a point $\tilde{E}$ if and only if the Zariski tangent space $T_{\tilde{E}}C$ is of dimension 1. If $E:=
\pi(\tilde{E})$, then the diagram \ref{cartesiano} yields a diagram
\begin{eqnarray} 
\shorthandoff{;:!?}
\xymatrix  @!0 @R=1.5cm @C=2.5cm  {
        T_{\tilde{E}}C \ar@{^{(}->}[r] \ar[d]_{d\pi}   &  T_{\tilde{E}}\tilde{X}^{(2)}  \ar[d]_{df^{(2)}}\\
        T_E\mathbb{P}^1    \ar@{^{(}->}[r]^{dj} & T_E X^{(2)}      
     }
\end{eqnarray} 
Since $\tilde{X}^{(2)}$ and $X^{(2)}$  are smooth of dimension 2, hence
\begin{eqnarray*}
\Dim T_{\tilde{E}}C &=& \Dim \Ker df^{(2)} + \Dim (\im dj \cap \im df^{(2)} ) \\
 &=& 2 - \Dim df^{(2)} + \Dim (\im dj \cap \im df^{(2)} ) \\
 &=& 2 +  \Dim \im dj  -  \Dim (\im dj + \im df^{(2)} ) \\
 &=& 1 + (2  -  \Dim (\im dj + \im df^{(2)} )). \\
\end{eqnarray*}
Then $C$ is smooth in $\tilde{E}$ if and only if $\im dj + \im df^{(2)}= T_E X^{(2)}$ and this happens if and only if the composition
\begin{displaymath}
 T_E\mathbb{P}^1  \stackrel{ dj}{\hookrightarrow} T_E X^{(2)} \rightarrow \mathop{\rm Coker}\nolimits  df^{(2)}
\end{displaymath}
is surjective. Recall the canonical isomorphisms (cf.\cite{ACGH} Ex.IV B-2) 
$$
\begin{array}{l}
T_E X^{(2)} \simeq H^0(X,\mathcal{O}_E (E)), \\ 
T_{\tilde{E}}\tilde{X}^{(2)} \simeq  H^0(\tilde{X},\mathcal{O}_{\tilde{E}} (\tilde{E})), \\
 T_E\mathbb{P}^1 \simeq H^0(X,\mathcal{O}_X (E)) /H^0(X,\mathcal{O}_X).
\end{array}
$$
The sheaf $\mathcal{O}_X(E)$ may be considered as a subsheaf of the sheaf of rational functions on $X$ by setting  $\mathcal{O}_X(E)_q = m^{-\nu_q(E)}_{X,q}$
for all $q \in X$, where $m^{-\nu_q(E)}_{X,q}$ is the ideal of rational functions with at most poles of degree $\nu_q(E)$ at $q$. Then  
$H^0(X,\mathcal{O}_X(E))$ is a subspace of the function field $K(X)$  and   $\mathcal{O}_E(E)_q = m^{-\nu_q(E)}_{X,q}/ \mathcal{O}_{X,q}$.
\\
Suppose $E = \sum_i n_iq_i$ with $q_i$ distinct points in $X$ and $n_i$ positive integers. For $h \in  H^0(\mathcal{O}_X (E))$ the $i$-th component
of $dj (h + H^0(\mathcal{O}_X))$ in  $\bigoplus m^{-n_i}_{X,q_i} =  H^0(\mathcal{O}_E (E)) $, is the image of $h$ in  $ m^{-n_i}_{X,q_i}/ \mathcal{O}_{X,q_i}$.
Denote by $q^1, q^2, q^3 \in \tilde{X}$ the points over $q$. Then, $\tilde{E} = \sum_i n_i^1q_i^1 + \sum_i n_i^2q_i^2 +\sum_i n_i^3q_i^3 $, with $n_i^1 +n_i^2 + 
n_i^3 = n_i$ and we can assume  $n_i^1 \geq n_i^2 \geq n_i^3$. Hence,
\begin{displaymath}
H^0(\mathcal{O}_{\tilde{E}} (\tilde{E})) = \bigoplus_i \bigoplus_{j=1}^3 m^{-n_i^j}_{X,q_i^j}/ \mathcal{O}_{X,q_i^j}.
\end{displaymath} 
Let $t_i$ be a local parameters of $X$ in $q_i$ and $t_i^j$ local parameters of $\tilde{X}$ in $q_i^j$ for $j=1,2,3$. Then $df^{(2)}((t_i^j)^{\nu}) =
(t_i)^{\nu}$ for $j=1,2,3$ and for all $\nu \in \mathbb{Z}$. Thus, $C$ is smooth in $\tilde{E}$ if and only if the map 
\begin{equation} \label{map}
H^0(\mathcal{O}_X (E)) \rightarrow \bigoplus_i  m^{-n_i}_{X,q_i}/   m^{-n_i^1}_{X,q_i}
\end{equation} 
is surjective. If $E =\pi(\tilde{E})$ is smooth or $E =\pi(2q_1)=2q$ , the right hand of \ref{map} is zero. If $\tilde{E}=q_1+q_2$, then $\pi(\tilde{E})=2q$ 
and the right hand of \ref{map} is $m^{-2}_{X,q}/   m^{-1}_{X,q}$. Let $h \in H^0(\mathcal{O}_X (E))$ with corresponding divisor $E_h \neq E$ in $g_2^1$. Then
$E_h = \Div (h) +E$, hence $\nu_q(E_h)=0=2+\nu_q(h)$, i.e., $\nu_q(h) = -2$. Therefore, the image of $h$ is a generator of  $m^{-2}_{X,q}/   m^{-1}_{X,q}$ and
$C$ is smooth in $E$.
\end{Dem}
{\bf Coverings of $\mathbb{P}^1$. }Let $X$ be a smooth projective curve of genus $g_X$ and consider a covering $f: X \rightarrow \mathbb{P}^1$ of degree
$n+2$, $n\geq 2$. Let $g_{n+2}^1$ be the linear system defining the covering and $B_f \subset \mathbb{P}^1$ the branch points. Assume that the  
$g_{n+2}^1$ is complete, then $n\leq g_X-1$. The degree of the ramification divisor, given by the Riemann--Hurwitz formula, is 
\begin{displaymath} 
w_f = 2g_X + 2n + 2.
\end{displaymath}
Let us to define set-theoretically the curve 
\begin{displaymath} 
C: = \{ E \in X^{(n)} \mid  |g_{n+2}^1 -E| \neq  \emptyset \}.
\end{displaymath}
For $z\in \mathbb{P}^1$, we denote $f^{-1}(z) = \{P_1 , \ldots,  P_{n+2}\}$ and by $P_{i_1, \ldots ,i_n} = P_{i_1} + \cdots +P_{i_n}$  a point in $C$, 
with $\{i_1, \ldots , i_n\} \subset \{ 1, \ldots , n+2 \}$. The curve $C$ comes with a map $h: C \rightarrow \mathbb{P}^1$ given by 
\begin{displaymath} 
P_{i_1,\ldots ,i_n} \mapsto f(P_{i_1}) = \cdots =  f(P_{i_n}) ,
\end{displaymath}
which is of degree ${n+2 \choose n} = \frac{(n+2)(n+1)}{2}$. In fact, $h$ corresponds to the composition of the classifying map $\pi_1(X-B_f, z_0) 
\rightarrow  \mathcal{S}_{n+2} $ with the monomorphism  $\mathcal{S}_{n+2} \hookrightarrow \mathcal{S}_N$ (with $N= \frac{(n+2)(n+1)}{2}$), which is given 
by the action of $\mathcal{S}_{n+2} $ on the cosets of the subgroup $\mathcal{S}_n \times \mathcal{S}_2 $, which  has a transitive image.  The proof of 
the fact that $C$ is smooth and irreducible is a slight 
generalization of Lemma 12.7.1 (\cite{B-L}).
\\
Let  $z \in  \mathbb{P}^1$ be a branch point above which the ramification is  $P_1 = P_2$, 
then the ramification of $h$ above $z$ is given by 
\begin{displaymath} 
P_{1,i_1,\ldots ,i_{n-1}} = P_{2,i_1, \ldots , i_{n-1}} \qquad \{ i_1, \ldots , i_{n-1}\} \subset \{3, \ldots , i_{n+2}\}.
\end{displaymath}
Hence the fibers of $f$ with only one ramification point induce ${n \choose n-1} = n$ simple ramification points for $h$. If above $z$ the covering $f$ has more 
than one simple ramification, the ramification  of $h$ above $z$ is a more complicated but not difficult to compute for  special cases. 
Using this information we can compute the degree of the ramification divisor $w_h$ of $h$ and the genus of $C$. For example, if the 
fibers of $f$ have no more than one simple ramification, we obtain 
\begin{equation} \label{simplecover} 
w_h = n\cdot w_f, \qquad g_C = n\cdot g_X + \frac{n(n+1)}{2}.
\end{equation} 
Let us to define an  effective symmetric correspondence on $C$ as follows
\begin{equation} \label{correspondencia}
D: P_{i_1,\ldots ,i_n} \mapsto  \sum P_{j_1,\ldots, j_n},
\end{equation} 
where the sum is over the subsets of $\{1, \ldots, n\}$ such that $|\{i_1, \ldots,i_n\}\cap \{j_1,\ldots ,j_n\}| = n-2 $. 
This is a correspondence of bidegree $(d,d)=(\frac{n(n-1)}{2},\frac{n(n-1)}{2} )$.
In order to $P_{1,\ldots,n}$ be a fixed point we must have 
\begin{displaymath} 
P_{1,\ldots,n} = P_{i_1,\ldots,i_{n-2}, n+1, n+2},
\end{displaymath}
for some $\{i_1,\ldots,i_{n-2} \} \subset \{1, \ldots , n\} $, then $P_{n+1} = P_{j_1}$ and  $P_{n+2} = P_{j_2}$ with $\{j_1, j_2\} \subset \{ 1, \ldots ,n\} $.
Then the correspondence has no fixed points on the fiber $h^{-1}(z) $ if and only if  $f^{-1}(z) $ admits at most one ramification point of index $\leq 3$.
In order to apply our criterion to this correspondence we must verify that 
the endomorphism induced by $D$, denoted by $\gamma_D$, satisfies the quadratic equation \ref{equation}.
\\
Put $D_j := D(P_{j_1, \ldots , j_{n-2},n+1,n+2})$,
where $j=\{j_1, \ldots , j_{n-2}\} \subset \{ 1, \ldots ,n\}$. We want to compute
\begin{displaymath} 
 D^2(P_{1,\ldots ,n})  =  \sum_j D_j.
\end{displaymath}
Observe that $P_{1,\ldots ,n}$ appears in $D_j$ for any $j$, then $P_{1,\ldots ,n}$ appears in  $D^2(P_{1,\ldots ,n})$ with coefficient $\frac{n(n-1)}{2}$.
Let us divide the elements of the fiber $h^{-1}(z)$ in two types, which are of the form $P_{k_1, \ldots , k_{n-2},n+1,n+2}$ and  
$P_{k_1, \ldots , k_{n-1},n+1}$ (or $P_{k_1, \ldots , k_{n-1},n+2}$ ).
Observe that $P_{k_1, \ldots , k_{n-2},n+1,n+2}$  appears in $D_j$ if and only if 
\begin{displaymath} 
|\{ k_1, \ldots , k_{n-2}\} \cap \{ j_1, \ldots , j_{n-2} \}| =n-4,
\end{displaymath} 
hence it is in as many $D_j$'s as subsets of $n-4$ elements of  $\{k_1, \ldots , k_{n-2}\}$, that is, a point of the form  $P_{k_1, \ldots , k_{n-2},n+1,n+2}$ 
appears in  $D^2(P_{1,\ldots ,n})$ with coefficient ${n-2 \choose n-4} = \frac{(n-2)(n-3)}{2}$.
\\
Similarly, $P_{k_1, \ldots , k_{n-1},n+1}$ appears in  $D_j$ if and only if 
\begin{displaymath} 
|\{ k_1, \ldots , k_{n-1}\} \cap \{ j_1, \ldots , j_{n-2} \}| =n-3,
\end{displaymath} 
hence $P_{k_1, \ldots , k_{n-1},n+1}$ appears in $D^2(P_{1,\ldots ,n})$ with coefficient ${n-1 \choose n-3} = \frac{(n-1)(n-2)}{2}$.
We can write
\begin{displaymath}
D^2(P_{1,\ldots ,n}) = \frac{n(n-1)}{2}P_{1,\ldots ,n} 
+ \frac{(n-2)(n-3)}{2} \sum_{\{ j_1 , \ldots, j_{n-2}\} \subset \{1, \ldots, n\}} 
P_{j_1, \ldots , j_{n-2},n+1,n+2} 
\end{displaymath}
\begin{displaymath}
 +  \frac{(n-1)(n-2)}{2} \left( \sum_{\{ j_1 , \ldots, j_{n-1}\} \subset \{1, \ldots, n\}} P_{j_1, \ldots , j_{n-1},n+1}
  + \sum_{\{ j_1 , \ldots, j_{n-1}\} \subset \{1, \ldots, n\}} P_{j_1, \ldots , j_{n-1},n+2} \right).  
\end{displaymath} 

Since  $  \frac{(n-1)(n-2)}{2} -\frac{(n-2)(n-3)}{2} = n-2 $ and  $  \frac{n(n-1)}{2} -\frac{(n-1)(n-2)}{2} = n-1 $ we have
\begin{displaymath}
D^2(P_{1,\ldots ,n}) = (n-1)P_{1,\ldots ,n} - (n-2)D(P_{1,\ldots ,n})  +  \frac{(n-1)(n-2)}{2}h^*(h(P_{1,\ldots ,n})).
\end{displaymath}
If $\gamma_D$ denotes the endomorphism of $JC$ induced by $D$, then we obtain the following equation in the Jacobian  
\begin{displaymath} 
\gamma_{D}^2 + (n-2)\gamma_D - (n-1) =0. 
\end{displaymath}
We have defined an effective symmetric correspondence on $C$ which verify the equation \ref{equation}.
\\
 Suppose that $D$ is free of fixed points. According to the Kanev's Criterion, the abelian variety $P:=\im (1-\gamma_D)$ is a Prym--Tyurin variety of exponent
 $n$ for the curve $C$. Using the formula \ref{dimension} and \ref{simplecover} to compute its dimension we obtain  $\dim P = g_X$. 
In fact, it was proved by Kanev (Proposition 8.5.4. \cite{Kanev2}) that $P$ is isomorphic to $JX$ as principal polarized abelian varieties.
We shall analyze some cases where the correspondence have fixed points.
\\
\\
\emph{Case n=2}. Suppose that we have a covering $f:X \rightarrow \mathbb{P}^1$ of degree 4 with two simple ramification points on two fibers and no more 
than one simple ramification point on the others fibers. The covering $h:C \rightarrow \mathbb{P}^1$ is of degree 6. In this case we have a correspondence on 
$C \times C$ of bidegree $(1,1)$, 
that is, an involution on the curve $C$, which sends $P_{1,2}$ to $P_{3,4}$. We get a fixed point of the correspondence for every fiber with two simple 
ramification points, then the involution has two fixed points. Computing the degree of ramification divisor of $w_f$ and $w_h$ we obtain the genus of 
$C$ as follows
\begin{eqnarray*} 
w_f=2g_X + 6\\
w_h = 2(w_f - 4) + 2(3)\\
g_C= -6 + \frac{w_h}{2} +1 =  2g_X .     
\end{eqnarray*} 
Applying the Theorem \ref{esencial} we have that $P:= \im (1-\gamma_D) \subset JC$  is a Prym--Tyurin variety of exponent 2. We use \ref{dimension} to compute its 
dimension 
\begin{eqnarray*} 
\dim P &=& \frac{1}{2}(g_C - 1 +1) \\
&=&  g_X.
\end{eqnarray*} 
\\
\emph{Case n=3}. Let   $f:X \rightarrow \mathbb{P}^1$ be a covering of degree 5 with ramifications as in the case $n$=2. The associated covering  
$h:C \rightarrow \mathbb{P}^1$ is 
of degree 10 and the correspondence $D$ on $C$ is of bidegree $(3,3)$. The fibers of $f$ with two simple ramification points, contribute with 3 
ramification points the index 4,2,2 on the corresponding fiber of $h$, indeed, if $P_1=P_2$ and $P_3=P_4$ are the ramification points on one of these 
fibers, we have on $C$
\begin{displaymath} 
P_{135}= P_{145}= P_{235}= P_{245}, \qquad P_{123}=P_{124}, \qquad P_{134}=P_{234}.
\end{displaymath}
Since $D(P_{135}) = P_{245}+P_{124}+P_{234} $, $ P_{135}$ is the only fixed point of $D$ on this fiber. Then $D$ has two fixed points. 
We have that
\begin{eqnarray*} 
w_f= 2g_X + 8\\
w_h = 3(w_f - 4) + 2(5)\\
g_C= -6 + \frac{w_h}{2} +1 =  3g_X + 2.     
\end{eqnarray*} 
By the Theorem \ref{esencial} we have that $P:= \im (1-\gamma_D)$  is a Prym--Tyurin variety of exponent 3. As before, we compute the dimension of $P$ and we obtain
 $\dim P = \frac{1}{3}(3g_X +2 -3 + 1) = g_X $.
\\
\\
\emph{Case n=4}. Consider a covering  $f:X \rightarrow \mathbb{P}^1$ of degree 6 having two fibers with three simple ramifications on each one and the others fibers 
with no  more than one simple ramification point. The associated covering $h:C \rightarrow \mathbb{P}^1$ is of degree 15 and the correspondence $D$ on $C$ is of bidegree $(6,6)$.
Let   $P_1=P_2$, $P_3=P_4$,  $P_5=P_6$ be the three ramification points on one fiber of $f$. They contribute, on the corresponding fiber of $h$, with
3 ramification points of index 4 as follows
\begin{eqnarray*} 
P_1 = P_{1235}  = P_{1245} = P_{1236} = P_{1246} \\
P_2 = P_{1345}  = P_{1346} = P_{2345} = P_{2346} \\
P_3 = P_{1356}  = P_{1456} = P_{2356} = P_{2456}. 
\end{eqnarray*} 
Observe that $P_1$,  $P_2$, $P_3$ are the three fixed points of $D$ on the fiber and they verify the conditions
b) of the Theorem \ref{esencial}. Applying the criterion we get a Prym--Tyurin variety $P$ of exponent 4. Observe that 
\begin{eqnarray*} 
w_f= 2g_X + 10\\
w_h = 4(w_f - 6) + 2(9)\\
g_C= 4g_X + 5.     
\end{eqnarray*}   
Hence the variety $P$ is of dimension $g_X$.
\\
\\
We shall show that in these examples $P \simeq JX$ as principally polarized abelian varieties. Let $Q_0 \in X$ and let $\alpha_n : X^{(n)} 
\rightarrow JX$ be the map
\begin{displaymath} 
E \mapsto \mathcal{O}_X (E-nQ_0), 
\end{displaymath}
for all $E \in X^{(n)} $. Let $\varphi = \alpha_{n_{|_C}} $ be the restriction to $C$.
\begin{Lem} \label{constante}
There exists a constant $b \in JX$ such that for all $Q \in C$ it verifies
\begin{displaymath} 
\varphi (D (Q))= (1-n) \varphi(Q) + b. 
\end{displaymath}
\end{Lem}
\begin{Dem}
Let
\begin{displaymath} 
b = \frac{n(n-1)}{2}(f^*(t) -(n+2)Q_0), 
\end{displaymath}
with $t \in \mathbb{P}^1$, then $b$ is a constant in $JX$ for any $t$. A straightforward computation shows that 
\begin{displaymath} 
\varphi (D (P_{1\ldots n})) + (n-1) \varphi (P_{1\ldots n}) = \frac{1}{2}n(n-1)(P_{1 \ldots n+2} - (n+2)Q_0). 
\end{displaymath}
 \end{Dem}
By the Universal Property of the Jacobian there exist a unique map $\tilde{\varphi}$ such that for all $c \in C$ the following diagram commutes
\begin{eqnarray}
\shorthandoff{;:!?}
\xymatrix @!0 @R=1.5cm @C=1.5cm  {
        C  \ar[d]_{\alpha_c} \ar[r]^{\varphi}  & JX \ar[d]^{t_{-\varphi (c)}} \\
	JC \ar[r]^{\tilde{\varphi}}  & JX     
     }
\end{eqnarray} 
The Lemma \ref{constante} tell us that the map $\tilde{\varphi}$ factorize by $P \subset JC$ since $\tilde{\varphi} \circ (\gamma_D +n - 1) =0 $. Then 
there exists a map $\psi$ making the following diagram commutative
\begin{eqnarray}
\shorthandoff{;:!?}
\xymatrix @!0  @R=1.5cm @C=1.5cm {
        C  \ar[d]_{\pi_c} \ar[r]^{\varphi}  & JX \ar[d]^{t_{-\varphi (c)}} \\
	P \ar[r]^{\psi}  & JX     
     }
\end{eqnarray} 
where $\pi_c = (1-\gamma_D) \circ \alpha_c$ is the Abel--Prym map.
\begin{Pro} \label{isomorfismo}
Let $P = \im (1 -\gamma_D) \subset JC$ with $\gamma_D$ the endomorphism induced by \ref{correspondencia}. If $\Dim P = \Dim JX$ then $\psi$ is an 
isomorphism of polarized abelian varieties.  
\end{Pro}
\begin{Dem}
Let $g= g_X= \Dim P$. Observe that $\varphi (C)$ generate $JX$ as abelian variety and since $Z$ and $JX$  have the same dimension $\psi : P 
\rightarrow JX$ is an isogeny. According to Welters' Criterion (Theorem 12.2.2 \cite{B-L}) we have
\begin{displaymath} 
\pi_*[C] = \frac{n}{(g-1)!} \wedge^{g-1} [\Xi]  \qquad \textrm{in} \qquad H^{2g-2}(P, \mathbb{Z}).  
\end{displaymath}
Suppose that have proved that
\begin{equation} \label {clase} 
\varphi_*[C] = \frac{n}{(g-1)!} \wedge^{g-1} [\Theta]  \qquad \textrm{in} \qquad H^{2g-2}(JX, \mathbb{Z}), 
\end{equation}
then $ \psi_* \wedge^{g-1} [\Xi] =  \wedge^{g-1} [\Theta]$ and by the Lemma 12.2.3. \cite{B-L} $\psi$ is an isomorphism.\\
Consider the sum map 
\begin{displaymath} 
s: X^{(2)} \times X^{(n)} \rightarrow X^{(n+2)}, 
\end{displaymath}
and define $Z:= s^{-1}(g^1_{n+2})$. Let $C' = \pi_1 (Z) \subset X^{(2)} $ and  $C = \pi_2 (Z) \subset X^{(n)} $ be the projections of $Z$ over the factors. 
The curves $C$ and $C'$ are isomorphic because both are isomorphic to $Z$. Then we  can describe $C'$ as follows   
\begin{displaymath} 
C' = \{ p+q \in  X^{(2)} \mid  |g^1_{n+2}-p-q| \neq \emptyset \}.
\end{displaymath}
\begin{Lem}
$\alpha_{2*}[C'] = \frac{n}{(g-1)!} \wedge^{(g-1)}[\Theta]$ in $H^{2g-2}(JX,\mathbb{Z})$.
\end{Lem}  
\begin {Dem}
It suffices to prove that $n[X]= \alpha_{2*}[C']$ in $H^{2g-2}(JX,\mathbb{Z})$ since the class of $X$ in $JX$ is $\frac{1}{(g-1)!} \wedge^{(g-1)}[\Theta]$. 
Let $q: X\times X \rightarrow X^{(2)}$ be the sum map and $\Delta_{\mathbb{P}^1}$ respectively $\Delta_X$ the diagonals in $\mathbb{P}^1$ respectively in $X$.
Since $[\Delta_{\mathbb{P}^1}] = [\mathbb{P}^1 \times \{ a\} ]+ [\{ b\} \times \mathbb{P}^1], \ a,b \in \mathbb{P}^1$, we have 
\begin{eqnarray*} 
[\Delta_X]  + q^*[C']& =  &(f \times f)^* ([\mathbb{P}^1 \times \{ a\} ]+ [\{ b\} \times \mathbb{P}^1]) \\
& = & [X \times D_a ]+ [D_ b \times X]\\
& =&  (n+2)[X \times \{ p \}] + (n+2)[\{ p\} \times X] \qquad \textrm{in} \quad H^2(X^2, \mathbb{Z}), 
\end{eqnarray*}
with $D_a, D_b \in g_{n+2}^1$ and for some  $p \in X$.
Let $\delta := q(\Delta _X)$. Applying $q_*$ we obtain
\begin{equation} \label{coh}
[\delta]  + 2[C'] = 2(n+2)[q(X \times \{ p \})]= 2(n+2)[X + p], 
\end{equation}
since $\Delta_X$ and $X \times \{ p\}$ are isomorphic to their images in $X^{(2)}$. Recall that $\alpha_{2*} (\delta) = 2_*[X] = 4[X]$ where $2_*$ is the push 
forward homomorphism of the multiplication by 2 in $JX$ (cf. Theorem 12.7.2 \cite{B-L}). Applying $\alpha_{2*}$ to the equation \ref{coh} we get
\begin{displaymath} 
4[X] + 2 \alpha_{2*}[C'] = 2(n+2) \alpha_{2*}[X + p]  =    2(n+2) [X ]  \qquad \textrm{in} \quad H^2(X^2, \mathbb{Z}). 
\end{displaymath}   
Thus,
\begin{displaymath} 
 \alpha_{2*}[C'] = n [X ]  \qquad \textrm{in} \quad H^2(X^2, \mathbb{Z}).
\end{displaymath}
This completes the proof.
\end{Dem}
We denote $\eta = \alpha_{n+2} (g^1_ {n+2}) \in W_{n+2} \subset JX$, where $W_{n} = \alpha_n(X^{(n)})$. Therefore, we have the diagram
\begin{eqnarray}
\shorthandoff{;:!?}
\xymatrix  @!0 @R=1.6cm @C=1.8cm  {
        & Z \subset X^{(2)} \times X^{(n)} \ar[rd]^{\pi_2} \ar[dl]_{\pi_1} \ar[rr]^s & &  X^{(n+2)} \supset g^1_{n+2} \ar[dd] \\
        C' \subset X^{(2)} \ar[d]^{\alpha_2}  & &  C \subset X^{(n)} \ar[d]^{\alpha_{n}} & \\
        \alpha_2(C') \subset W_2 \ar@{^{(}->}[rr] && W_n \ar@{^{(}->}[r] & W_{n+2} \ni \eta           
     }
\end{eqnarray} 
and then $(\eta - W_n).W_2 = \alpha_2(C')$ and $(\eta - W_2).W_n = \alpha_n(C)$. Hence
\begin{displaymath} 
\alpha_2(C') = \eta -\alpha_n(C), 
\end{displaymath} 
and we conclude that $\alpha_n(C)$ is algebraically equivalent to $n X$. This shows the equality ~\ref{clase} and the proof of the Proposition \ref{isomorfismo} 
is completed.
\end{Dem}
\begin{Rmk}
The examples coming from coverings of $\mathbb{P}^1$ are generalizations of the Recillas' construction \cite{R} 
\end{Rmk}

\end{document}